\documentclass[preprint,12pt]{elsarticle}
\usepackage{amssymb}
\usepackage{amsthm}
\usepackage{subcaption,verbatim}
\usepackage{amsfonts, amssymb, amscd}
\usepackage{mathtools, cuted}
\usepackage{hyperref}
\usepackage{multicol}
\usepackage{tikz}
\usepackage{ragged2e}

\journal{Number Theory}
\newtheorem{Def}{Definition}[section]
\newtheorem{Lem}{Lemma}[section]
\newtheorem{Thm}{Theorem}[section]
\newtheorem{Rem}{Remark}[section]
\newtheorem{Prop}{Proposition}[section]
\newtheorem{Ex}{Example}[section]
\newtheorem{Con}{Conjecture}[section]
\numberwithin{equation}{section}
\newtheorem{Cor}{Corollary}[section]
\begin{document}

\begin{frontmatter}

\title{The Mobius Function and Congruent Numbers}

\author{Roy Burson}

\begin{abstract}
This work provides a complete characterization of congruent numbers in terms of Pythagorean triples. Specifically, we show that every congruent number can be written as $$\frac{nm\left(m-n\right)\left(m+n\right)}{\sigma^2}$$ were as $$\sigma \vert \rho\biggl(\left(m-n\right)\left(m+n\right)\biggr),\indent \text{or} \indent \sigma \vert \rho( nm )$$ were $\rho(\alpha)$ denotes the non-square free part of its argument $\alpha$. As a consequence, in order to find congruent numbers it suffices to devise a condition so that the equality $\mu(m-n)+1 = \gcd(m,n)$ or $\mu(m+n)+1 =\gcd(m,n)$ holds, were $\mu$ is the Mobius function.
\end{abstract}

\begin{keyword}
Congruent Numbers \sep The Mobius Function 


\end{keyword}

\end{frontmatter}


\section{Introduction}
\label{S1}
According to \cite{Arab} the set of congruent numbers is one of the oldest sets of numbers on record first appearing in the Arab Manuscripts. The most recent results concerning congruent numbers stems from the properties of elliptic curves \cite{GI}, which is connected to the Birch and Swinnerton-Dyer Conjecture discussed by \cite{Wiles}. Assuming the validity of the Birch and Swinnerton-Dyer Conjecture conjecture \cite{Tunnel} derived a criteria to determine if a given number $n$ is congruent. Despite of this result, the congruent number problem has yet to be solved. The result in \cite{Tunnel} has brought a wide spread attention to the properties of elliptic curves. \\
\indent Due to this wide spread attention surrounding the group properties of elliptic curves no new elementary results have been seen in the recent years. Instead of developing any new ideas most authors have pushed there studies toward a solution of the Birch and Swinnerton-Dyer Conjecture by studying the group properties of Elliptic curves and attempting to address the rank, which is not the case in this paper. In this work we address the congruent problem by developing some new results, which have not been seen in other literature. The results obtained here are elementary as they are independent of any previous results surrounding congruent numbers except for Euclid's theorem that establishes a parametric description of Improper-Congruent numbers (which we discuss later) and \cite{Conrad}'s quick discuss on how to find square free congruent numbers.  \cite{Conrad} actually suggest that using the Pythagorean theorem is a bad idea to approach the congruent number problem, which this work proves wrong by showing that congruent numbers are specific divisors of Pythagorean triples.  \\
\indent This paper serves to provide the complete parametric description of all congruent numbers using Pythagorean triples. Specifically, this work shows that all congruent numbers are either Pythagorean triples or they are the non-square free divisors of Pythagorean triples. Hence, in order to find congruent number this work shows it suffices to examine the non-square free divisors of the Pythagorean triples. Thus providing a direct link between Pythagorean triples and all congruent numbers. Although this work provided the complete analytical description of congruent numbers non of these results are connected to Elliptic curves as it would overcompensate the procedure presented herein as it is long enough.\\
\indent The main objective of this work is too separate the congruent numbers into three different categories: Improper, Semi-Proper, and Proper. In each case we gather some results independently and then bring them all together towards the end of the work.   Improper congruent numbers are those congruent numbers which represent the area of a triangle which has both legs strictly in $\mathbb{Z}$. Semi-Proper congruent numbers are those congruent numbers which represent the area of a triangle with one leg in $\mathbb{Z}$ and the other strictly in $\mathbb{Q}$. Proper-Congruent numbers are those congruent numbers which has all sides in $\mathbb{Q}$ and strictly not in $\mathbb{Z}$.\\
\indent In section 2 we develop some preliminary results that are needed throughout the work. Moreover the ideas in this section have shaped or formed the rest of the work. In this section we discuss Euclid's famous theorem (Theorem $\S$ \ref{Thm 2.1}) regarding the parametric description of Pythagorean triples and show how it can be used to develop some results concerning proper congruent numbers. The main result in this section is provided by Theorem $\S$ \ref{Thm 2.2}, which provides a condition to find the Proper-Congruent numbers. This result combined with Proposition $\S$ \ref{Prop 2.2} shows that there is a correspondence between the proper congruent numbers and the divisors of the Pythagorean triple, which is discussed more thoroughly in Section 7. \\
\indent Section 3 shows how one can use the Mobius function to generate some congruent numbers. This section established some preliminaries which help establish a connection between the proper congruent numbers and the altered Fermat equation $ax^4+by^4=z^2$ discusses in the latter section.
\section{Primitive Triples and Proper Congruent Numbers}
\label{S2}
\begin{Def}
\label{Def 2.1}
A \textit{Pythagorean triple} is a set of integers $\{a,b,c\}$ so that   $a^2+b^2 = c^2$ were as $a,b,c\in \mathbb{Z}$ (for, example $(3,4,5)$). If $\gcd(a,b,c) = 1$ we say that $\{a,b,c\}$ is a Primitive triple.
\end{Def}
\begin{Thm}(Euclid's Elements \cite{RFEE}). 
\label{Thm 2.1}
All primitive Pythagorean triples $(a,b,c)$ have the representation \begin{align*}
a & = m^2-n^2\\
b& = 2nm\\
c &= m^2+n^2
\end{align*}
were $\gcd(n,m)=1$ and $m>n$ and $m\not \equiv n(\bmod{2})$.
\end{Thm}
\begin{Def}
\label{Def 2.2}
A \textit{congruent number} is a number $n\in \mathbb{N}$ that is the area of some right triangle $\triangle ABC$ with rational sides $(\frac{a}{b},\frac{c}{d},\frac{e}{f})$.
\end{Def}
\begin{Def}
\label{Def 2.3}
A \textit{Proper-Congruent number} is a number $n\in \mathbb{N}$ that represents the area of some right triangle $\triangle ABC$ with rational sides $(\frac{a}{b},\frac{c}{d},\frac{e}{f})$ and the additional property that neither $b$ nor $d$ can be one. That is $b\neq 1 \neq d$. So both sides lengths of $\triangle ABC$ are strictly in $\mathbb{Q}$ and not $\mathbb{Z}$.
\end{Def}
\begin{Prop}
\label{Prop 2.1}
If $n$ is a Proper-Congruent number with $\triangle ABC$ constructed by the triple $(\frac{a}{b},\frac{c}{d},\frac{e}{f})$, then $f=bd$. \begin{proof}
Suppose $n$ is a proper congruent number with $1=\gcd(a,b)=\gcd(c,b)=\gcd(e,f)$. Since $n$ is congruent it represents the area of the triangle $\triangle ABC$ which has area $n=\frac{1}{2}(\frac{a}{b})(\frac{c}{d}) = \frac{ac}{bd}$. Now $n$ is an integer so $bd\vert ac$. Since $n$ is proper we know that $b$ does not divide $c$ and $d$ does not divide $a$. Therefore $d\vert a$ and $b\vert c$. Let $i$ denote the largest power of $b$ so that $b^{i}\vert c$ but $b^{i+1} \not \vert c$. Similarly, let $j$ denote the largest power of $d$ so that $d^{j}\vert a$ but $d^{j+1} \not \vert a$. Write $a =d^{j}q_2$ and $c =b^iq_1$ for some quotients $q_1$ and $q_2$, then by the hypothesis it follows that 
\begin{align*}
(\frac{a}{b})^2+(\frac{c}{d})^2 & =(\frac{e}{f})^2\\
& \Downarrow\\
(\frac{d^{j}q_2}{b})^2+(\frac{b^iq_1}{d})^2 & =(\frac{e}{f})^2\\
& \Downarrow\\
\frac{d^{2j}q_2^2}{b^2}+\frac{b^{2i}q_1^2}{d^2} & =\frac{e^2}{f^2}\\
& \Downarrow\\
\frac{d^{2j+1}q_2^2+ b^{2i+2}q_1^2}{b^2d^2} & =\frac{e^2}{f^2}\\
\end{align*}
Now suppose that the numerator and the denominator of the last equality are not equal. Then it must follow that $$b^2d^2\vert \left(d^{2j+1}q_2^2+ b^{2i+2}q_1^2\right)$$ but $bd\vert b^2d^2$ so then $$bd\vert b^2d^2\vert  \left(d^{2j+1}q_2^2+ b^{2i+2}q_1^2\right)$$ Thus product $bd$ divides the numerator. This implies $b$ and $d$ divide the numerator individually. This clearly cannot happen because $b$ divides the numerator if and only if $b\vert q_1^2$, but $a=d^jq_1$ and $1=\gcd(a,b)=\gcd(b,d^jq_1)$ so $b$ cannot divide $q_1$. The same argument holds for the integer $d$. Thus the numerator and the denominator are equal, which implies that is $f^2 = b^2d^2 = (bd)^2$. Hence, $f=bd$ as desired. 
\end{proof}
\end{Prop}

\begin{Prop}
\label{Prop 2.2}
If $n$ is a proper congruent number with $\triangle ABC$ constructed by the triple $(\frac{a}{b},\frac{c}{d},\frac{e}{f})$, then $e^2 = (da)^2+(bc)^2$. \begin{proof}
Assume that $n$ is a congruent number with the right triangle $\triangle ABC$ constructed by the triple $(\frac{a}{b},\frac{b}{c},\frac{e}{f})$. By proposition $\S$ \ref{Prop 2.1} it follows that $f=bd$. Therefore, \begin{align*}
(\frac{e}{f})^2 & = (\frac{a}{b})^2+(\frac{c}{d})^2 \\
& \Downarrow\\
(\frac{e}{bd})^2 & = (\frac{a}{b})^2+(\frac{c}{d})^2 \\
& \Downarrow\\
e^2  & = (bd)^2\left((\frac{a}{b})^2+(\frac{c}{d})^2\right) \\
& \Downarrow\\
e^2  & = (bd)^2(\frac{a}{b})^2+(bd)^2(\frac{c}{d})^2 \\
& \Downarrow\\
e^2  & = (da)^2+(bc)^2 
\end{align*}
\end{proof}
\end{Prop}
\begin{Rem}
\label{Rem 2.1}
One is tempted to believe that Proposition $\S$ \ref{Prop 2.1} and  $\S$ \ref{Prop 2.2} follow directly from the Pythagorean theorem. However, this is not true. For example, if we have a right triangle $\triangle ABC$ with sides $(\frac{a}{b},\frac{c}{d},\frac{e}{f})$ then the Pythagorean theorem says that $(\frac{e}{f})^2 = (\frac{a}{b})^2+(\frac{c}{d})^2 = \frac{(ad)^2+(bc)^2}{(bd)^2}$. One may now be attempted to assume that $e = \sqrt{(ad)^2+(bc)^2}$ and $f =bd$ but we cannot jump directly to this conclusion because this is not a general result that holds. For example, if $\frac{x}{y}=\frac{k}{l}$ then $x=k$ and $y=l$ if and only if $\gcd(x,y)=1=\gcd(k,l)$, and here we do not know the greatest common factor so we may not make this conclusion. However, it turns out to be true as we will demonstrate below.
\end{Rem}
\begin{Prop}
\label{Prop 2.3} If $(a,b,c)$ is a primitive Pythagorean triple then, $\gcd(a,b)=1$. \begin{proof} Suppose $(a,b,c)$ is a primitive triple with $\gcd(a,b) \neq 1$. Then set $d = \gcd(a,b) >1$. By the Euclidean division algorithm there exist $k_1$, $k_2$ so that $a = dk_1$ and $b = dk_2$. Then \begin{align*}
c^2 & = a^2 +b^2 \\
& \Downarrow\\
c^2 & = (dk_1)^2 +(dk_2)^2 \\
& \Downarrow\\
c^2 & = d^2(k_1^2 +k_2^2)\\
& \Downarrow\\
(\frac{c}{d})^2 & = k_1^2 +k_2^2
\end{align*}
Since $k_1$ and $k_2$ are integers it must follow that $d\vert c$. But this means $\gcd(a,b)\vert c$. Since $\gcd(a,b)\vert a$ and $\gcd(a,b)\vert b$ it follows that $d = \gcd(a,b)$ is a divisor or $a$, $b$, and $c$. So $\gcd(a,b,c)\ge \gcd(a,b) >1$. So then $\gcd(a,b,c)>1$ which is to say that $(a,b,c)$ is not primitive. This is a contradiction.
\end{proof}
\end{Prop}
\begin{Thm}
\label{Thm 2.2}
Let $a,b,c,d,e$ be integers that satisfy the following:\\
\noindent  (i) $e^2 = (da)^2+(bc)^2$\\
(ii) $d\vert a$, $b\vert c$, $\gcd(a,b) = 1=\gcd(c,d)$, and $b\neq 1\neq d$.\\
(iii) $\{a,c\}$ contains an even integer,  and this even number has more factors of 2 than $b$ or $d$.\\
Then the number $\frac{1}{2}(\frac{a}{b})(\frac{c}{d})$ is a proper congruent.
\begin{proof}
Assume $(i)-(iii)$ are valid. Set $f=bd$ and observe that \begin{align*}
e^2 & = (da)^2 +(bc)^2 \\
& \Downarrow\\
\frac{e^2}{f^2} & = \frac{(da)^2 +(bc)^2}{f^2} \\
& \Downarrow\\
\frac{e^2}{f^2} & = \frac{d^2a^2}{f^2}+\frac{b^2c^2}{f^2} \\
& \Downarrow\\
\frac{e^2}{f^2} & = \frac{a^2}{b^2}+\frac{c^2}{d^2} \\
& \Downarrow\\
(\frac{e}{f})^2 & = (\frac{a}{b})^2+(\frac{c}{d})^2 \\
\end{align*}
This is a rational right triangle with area $\frac{1}{2}(\frac{a}{b})(\frac{c}{d})$. This is a whole integer since $d\vert a$ and $b\vert c$ and by condition (iii) the number  $a$ or $c$ contains at least one more factor of $2$ its divisor $b$ or $d$, which also satisfy the condition $b\neq 1\neq d$. We only need to very that $\gcd(e,f)=1=\gcd(a,b)=\gcd(c,d)$. Condition (ii) gives that $1=\gcd(a,b)=\gcd(c,d)$ so we only need to show $\gcd(e,f)=1$. By condition (i) we can see that $bd$ does not divide $e^2=$ since $e^2=(da)^2+(bc)^2$. So $bd$ cannot divide $e$, which means $\gcd(f,e)=\gcd(bd,e) = 1$. 
\end{proof}
\end{Thm}
\begin{Ex}
\label{Ex 2.1}
Let us try to find a Pythagorean triple that satisfies the conditions of the hypothesis in Theorem $\S$ \ref{Thm 2.1} and $\S$ \ref{Thm 2.2}. According to \cite{alvaro}, around 1220 Fibonacci was challenged by Johannes Palermo to find a right triangle whose are is equal to $5$, and Fibonacci found $\{\frac{3}{2},\frac{20}{3},\frac{41}{6}\}$. Nobody  actually knows how he came up with this triangle, however, he stated that if $n$ is perfect square then $n$ is not a congruent number. The proof had to wait until (1601-1665) when Pierre de Fermat finally solved the problem.\\
\indent Now let us attempt to find how Fibonacci found this triangle using Theorem $\S$ \ref{Thm 2.1} and $\S$ \ref{Thm 2.2}. By Theorem $\S$ \ref{Thm 2.1} it follows that all the primitive Pythagorean triples $(a,b,c)$ have the representation
\begin{align*}
a & = m^2-n^2\\
b& = 2nm\\
c &= m^2+n^2
\end{align*}
were $\gcd(n,m)=1$ and $m>n$ and $m\not \equiv n(\bmod{2})$. Choose $m=5,n=4$ then $\gcd(5,4)=1$, $5\neq (4{\bmod{2}})$, and $5>4$ so these integers qualify to an integer Pythagorean triple. That is $a = (m-n)(m+n)=1\cdot 9=9$, $b = 2\cdot m\cdot n= 2\cdot 3\cdot 2 = 12$. Both $a$ and $b$ are square free so we can break them up by factoring the largest non square free part of both $a$ and $b$ as $a=3\cdot 3=9$ and $b=2\cdot 20 $. Write $c^2 = (9)^2+(40)^2 = (3\cdot3)^2 +(2\cdot 20)^2$, then by direct application of Theorem $\S$ \ref{Thm 2.2} it follows that the integer $n = \frac{1}{2}(\frac{3}{2})(\frac{20}{3}) =5$ is a congruent number generated by the triple $(\frac{3}{2},\frac{20}{3}, \frac{41}{6})$. This is precisely the exact triple that Fibonacci discovered in 1220 in response to his challenger. 
\end{Ex}
\begin{Def}
\label{Def 2.4}
We say a number $n$ is square free if $p^2\not \vert n$ for all $p\vert n$. A number is not square free if there is a divisor $p\neq 1$ of $n$ so that $p^2\vert n$.
\end{Def}
\begin{Lem}
\label{Lem 2.1}
Assume that $m,n\in \mathbb{Z}$. If $\gcd(m,n)=1$, then $\gcd(m-n,m+n)=1$ or $\gcd(m-n,m+n)=2$ . \begin{proof}
Assume that $m,n\in \mathbb{Z}$ were as $\gcd(m,n)=1$. Let $d$ be a common positive divisor of $m+n$ and $m-n$. Then $d\vert m+n$ and $d\vert m-n$. Hence, $d\vert (m+n)+(m-n)=2m$. Similarly, $d\vert (m+n)-(m-n)=2n$. So we have that $d\vert  2n$ and $d\vert 2m$. Since $\gcd(m,n)=1$ $d$ must divide $2$. So then $d=1$ or $d=2$.
\end{proof}
\end{Lem}
\section{Generation}
\label{S3}
\begin{Def}
\label{Def 3.1}
Let $n\in \mathbb{N}$. The Mobius function is the map $$ \mu(n)=
\begin{cases}
1 & \textit{if $n$ is square free}\\
0 & \textit{if $n$ is not square free}\\
\end{cases}$$ 
\end{Def}
\begin{Con}
\label{Con 3.1}
If $n\in 2^{\mathbb{N}}$ then $n$ is not congruent.  
\end{Con}
\begin{Prop}
\label{Prop 3.1}
Let $p$ be a prime number. If there exist $(k,l)$ so that $k^2-64p^2 = l^2$ were as $l >1$ then $p$ is a congruent number. \begin{proof}
Suppose $p$ is a prime number and suppose that there exist $(k,l)$ so that $k^2-64p^2 = l^2$. Write $k^2=l\cdot l+(2\cdot (4p))^2$. Define $d=l=a$,$b=2$, and $c=4p$. Since $\gcd(l,8p)=1$ it follows that $\gcd(d,4p)=1$. Since $d=l>1$ we see that $d$ does not divide $p$ (otherwise we would have a contradiction). Thus $\gcd(a,b)=1=\gcd(c,d)$. Now by direct application of Theorem $\S$ \ref{Thm 2.2} it follows that the number $n = \frac{1}{2}(\frac{a}{d})(\frac{4p}{2})$ is congruent. But since $a=d$ then we trivially have $n=p$. Hence, we have shown that $p$ is congruent.
\end{proof}
\end{Prop}
\begin{Prop}
\label{prop 3.2}
If $|m-n|=1$, then $\gcd(m-n,m+n)=1$ \begin{proof}
Let $|m-n|=1$ so that $m=n+1$ or $n=m+1$. In the first case we have $\gcd(m-n,m+n)=\gcd((n+1)-n,(m+1)+m) = \gcd
(1,2m+1)=1$. Similarly, in the second case, if $n=m+1$ then, $\gcd(m-n,m+n)=\gcd(m-(m+1),n+(n-1)) = \gcd(1,2n-1)=1.
(1,2n+1)$. These are the only possibilities so we are done.
\end{proof}
\end{Prop}
\begin{Cor}
\label{Cor 3.1}
If $\gcd(m-n,m+n)=2$, then $|m-n|\neq 1$ \begin{proof}
This is the contropositive of proposition $\S$ \ref{prop 3.2}. By Lemma $\S$ \ref{Lem 2.1} if $\gcd(m-n,m+n)\neq 1$ it must be that $\gcd(m-n,m+n) = 2$. Since this is the only other possibility.
\end{proof}
\end{Cor}
\begin{Prop}\label{prop 3.3}
If $\gcd(m,n)=1$ and $m\not \equiv n(\bmod{2})$ then $\gcd(m-n,m+n)=1$. \begin{proof}
Assume $\gcd(m,n)=1$ and $m\not \equiv n(\bmod{2})$. By lemma $\S$ \ref{Lem 2.1} it follows that $\gcd(m,n)=1$ or $\gcd(m,n)=2$. Since $m\not \equiv n(\bmod{2})$ it follows that $m-n$ and $m+n$ are both odd. Hence, $2$ cannot divide neither of these numbers. This leaves the only conclusion that $\gcd(m,n)=1$. 
\end{proof}
\end{Prop}
\begin{Rem}
\label{Rem 3.1}
Theorem $\S$ \ref{Thm 2.2} illustrates that there is some type of correspondence between non-square free Pythagorean triples and congruent numbers (not necessarily a bijection). In this sense, one can generate all some congruent numbers by finding all non-square free Pythagorean triples. Later in definition $\S$ \ref{Def 2.1} we classify these types of congruent as \textit{Proper Congruent Numbers}. Recall Euclid's Theorem \ref{Thm 2.1} asserts that all Pythagorean triples  $\{a,b,c\}$ have the form
 \begin{align*}
a & = m^2-n^2\\
b& = 2nm\\
c &= m^2+n^2
\end{align*}
were $\gcd(n,m)=1$ and $m>n$ and $m\not \equiv n(\bmod{2})$. Since $m\not \equiv n(\bmod{2})$ it follows that $4\vert b$ and hence $\mu(b)=0$. This indicates that one needs only to find when $a=m^2-n^2$ is not square free. But we now $a=m^2-n^2=(m-n)(m+n)$ so that $\mu(a)=0$ whenever $\mu(m-n)=0$, $\mu(m+n)=0$, or $\mu(m^2-n^2)=0$. In the next proposition we will show that $\mu(m^2-n^2)=0$ can not happen,which leaves us to find when $m-n$ or $m+n$ is not square free.
\end{Rem}
\begin{Ex}
\label{Example 3.1}
Consider $m=2^2\cdot5^2\cdot7^3=34300$ and $n=11^2\cdot 13^2=20449$. Then set $\beta = 2mn=(2^3\cdot5^2\cdot 7^3)(11^2\cdot 13^2)$ and compute $\alpha = m^2-n^2=(m-n)(m+n)=(13851)(54749)=(3^6\cdot 19)(53\cdot 1033)$. Clearly, $\mu(\beta)=0$ since $\mu(n)=0$ and $\mu(m)=0$, and $\mu(\alpha)=0$ since $3^6\vert \alpha$. By Euclid's Theorem $\S$ \ref{Thm 2.1} we now $(\alpha,\beta,\gamma)$ is a primitive Pythagorean triple. Henceforth, using Theorems $\S$ \ref{Thm 2.2} we can find a congruent number that corresponds to this triple (Note! the theorem does not state that the congruent number is unique. This is because on some occasions we may form more than one congruent number for a given primitive triple $\{\alpha,\beta,\gamma\}$ (see table 4.1).\\
\indent To find a corresponding congruent number we may take any divisors $a,d\vert \alpha$ and $b,c\vert \beta$ so that $(m^2+n^2)^2=(da)^2+(bc)^2$ were $d\vert a$, $b\vert c$. We can do this by rearranging the factorization of $\alpha$ and $\beta$. For example, set $d=3$, $a=3^5\cdot 19$, $b=2\cdot 5\cdot 7$, and $c=2^2\cdot 5\cdot 7^2\cdot 11^2\cdot 13^2$. Allow $e=m^2+n^2$ then $e^2=(da)^2+(bc)^2$ were $d\vert a$, $b\vert c$ and $d\neq 1\neq b$. Corollary to Theorem $\S$ \ref{Thm 2.2} it follows that the number $n= \frac{1}{2}(\frac{a}{b})(\frac{c}{d})=\frac{1}{2}(\frac{3^5\cdot 19}{2\cdot 5\cdot 7})(\frac{2^2\cdot 5\cdot 7^2\cdot 11^2\cdot 13^2}{3})= 2\cdot 3^4\cdot 7\cdot 11^2\cdot 13^2=23189166$ is congruent. We could have also chosen $d=3$, $a=3^5\cdot 19$, $b=2$, and $c=2^2\cdot 5^2\cdot 7^3\cdot 11^2\cdot 13^2$, so that we would have found $n = \frac{1}{2}(\frac{a}{b})(\frac{c}{d})=\frac{1}{2}(\frac{3^5\cdot 19}{2})(\frac{2^2\cdot 5^2\cdot 7^3\cdot 11^2\cdot 13^2}{3}) = 1092566475$ instead of $23189166$ as above. To find all congruent numbers that correspond to this specific Pythagorean triple we can list all the different factorization and compute each different congruent number. The results of each possible triple is illustrated in Table \ref{table 1} below.
\begin{table}[h!] 
\centering
\begin{tabular}{c|c}
$n$ & $\frac{1}{2}(\frac{a}{b})(\frac{c}{d})$\\
\hline
$1943020219000$ & $\frac{1}{2}(\frac{3^6\cdot 19}{1})(\frac{2^3\cdot 5^2\cdot 7^3\cdot 11^2\cdot 13^2}{1})$\\
$485755054800$ & $\frac{1}{2}(\frac{3^6\cdot 19}{2})(\frac{2^2\cdot 5^2\cdot 7^3\cdot 11^2\cdot 13^2}{1})$\\
$215891135500$ & $\frac{1}{2}(\frac{3^5\cdot 19}{1})(\frac{2^2\cdot 5^2\cdot 7^3\cdot 11^2\cdot 13^2}{3})$\\
$53972783870$ & $\frac{1}{2}(\frac{3^5\cdot 19}{2})(\frac{2^2\cdot 5^2\cdot 7^3\cdot 11^2\cdot 13^2}{3})$\\
$5996975985$ & $\frac{1}{2}(\frac{3^4\cdot 19}{2})(\frac{2^2\cdot 5^2\cdot 7^3\cdot 11^2\cdot 13^2}{3^2})$\\
$122387265$ & $\frac{1}{2}(\frac{3^3\cdot 19}{2})(\frac{2^2\cdot 5^2\cdot 7^3\cdot 11^2\cdot 13^2}{3^3})$\\
$1011465$ & $\frac{1}{2}(\frac{3^3\cdot 19}{2\cdot 5})(\frac{2^2\cdot 5\cdot 7^3\cdot 11^2\cdot 13^2}{3^3})$\\
$5985$ & $\frac{1}{2}(\frac{3^3\cdot 19}{2\cdot 5\cdot 7})(\frac{2^2\cdot 5\cdot 7^2\cdot 11^2\cdot 13^2}{3^3})$\\
$8645$ & $\frac{1}{2}(\frac{3^3\cdot 19}{2\cdot 5\cdot 7\cdot 11})(\frac{2^2\cdot 5\cdot 7^2\cdot 11\cdot 13^2}{3^3})$\\
$665$ & $\frac{1}{2}(\frac{3^3\cdot 19}{2\cdot 5\cdot 7\cdot 11\cdot 13})(\frac{2^2\cdot 5\cdot 7^2\cdot 11\cdot 13}{3^3})$\\
\end{tabular}
\caption*{Table 3.1. List of different congruent numbers that correspond to the triple $(758328399,200400200,615222200900000000)$}
\end{table}\label{table 1}
To see why this is useful one has to notice that with a single Pythagorean triple we were able to find $10$ congruent numbers shuffling the factorization of each term in $\frac{1}{2}(\frac{a}{b})(\frac{c}{d})$. Ideally, one would like to find large Pythagorean triples $(\alpha,\beta,\gamma)$ were $\mu(\alpha)=0=\mu(b)$ and the number of different factors of $a$ and $b$ is large. This is explored later in section 6 and 7.
\end{Ex}

\begin{Def}\label{Defintion 3.3}
Suppose that $n\in \mathbb{N}$ and $a\vert n$ were $a\neq 1$. The mapping $\varrho:\mathbb{N}\rightarrow \mathbb{N}$ defined by the rule $n\mapsto\max\{d: d\vert n, ~n= d^2a, ~a\neq 1 \}$ is called the \textit{Square Part} of $n$. It is the non-square free part of an integer. If $n$ is square free then we decree that $\varrho(n)=1$.
\end{Def}
\begin{Lem}\label{Lemma 3.1}
let $p$ be a prime number then $\varrho(p^i) = p^{\lfloor\frac{i}{2}\rfloor}$ for each $i\ge 2 \in \mathbb{N}$. \begin{proof}
let $p$ be a prime number and consider $y=p^i$ for $i\ge 2\in \mathbb{N}$. Now by definition $\varrho(p^i)=\max\{d: d\vert p^i, ~p^i= d^2a, ~a\neq 1 \}$. But any divisor of $p$ is either $1$ or $p$. So the maximum divisor $d$ that fits this property must by at least $p$ if $i\ge 2$. Note that $\left(p^{\lfloor\frac{i}{2}\rfloor}\right)\left(p^{\lfloor\frac{i}{2}\rfloor}\right) =p^{\lfloor\frac{i}{2}\rfloor+\lfloor\frac{i}{2}\rfloor}=p^{\lfloor\frac{2i}{2}\rfloor} =p^i$ and $\left(p^{\lfloor\frac{i}{2}\rfloor}\right)$ is the largest integer that divides $\left(p^{\lfloor\frac{i}{2}\rfloor}\right)$. Hence, $\varrho(p^i) = p^{\lfloor\frac{i}{2}\rfloor}$ and we are done.
\end{proof}
\end{Lem}
\begin{Lem}\label{Lemma 3.2}
The \textit{Square Part} of $n$ is multiplicative over relatively prime numbers. That is if $a,b\in \mathbb{N}$ so that $\gcd(a,b)=1$ then $\varrho(a\cdot b) = \varrho(a)\cdot \varrho(b)$. \begin{proof}
Let $\gcd(a,b)=1$. If $\varrho(a)=1$ or $\varrho(b)=1$ then the result is trivial. So assume that $\varrho(a)\neq 1$ and $\varrho(b)\neq 1$. This being the case we can write $a = A^2k$ and $b=B^2l$ were $\varrho(a)= A^2$ and $\varrho(b)=B^2$. Then $ab = A^2B^2kl=(AB)^2(kl)$. Moreover, since $A^2\vert k$ and $B^2\vert l$ it follows that $(AB)^2 \vert (kl)$. Since $\gcd(a,b)=1$ it must follow that $\gcd(AB,kl)=1$. Thus the square part of $ab$ is $(AB)^2$ and hence $\varrho(a)\cdot \varrho(b)=A^2B^2 = (AB)^2=\varrho(ab)$ as desired.
\end{proof} 
\end{Lem}
\begin{Prop}\label{Proposition 3.4} Let $n\in \mathbb{N}$ and write $n$ in conical form $n=\prod_{i=1}^{k}{p^{e_i}_i}$ for primes $p_1,p_2,\cdots,p_k$ and integers $e_1,e_2,\cdots,e_k\in \mathbb{N}$. Then $$\varrho(n)=\prod_{i=1}^{k}{p^{\lfloor\frac{e_i}{2}\rfloor}_i}$$ \begin{proof}
This proof utilities the multiplicative property of the $\varrho$. Write $n=\prod_{i=1}^{k}{p_i^{e_i}}$ it its conical factorization. Since $p_i$ is prime for each $i\in \{1,2,\cdots ,k\}$ it follows that $\gcd(p_i^{e_i},p_{j}^{e_j})=1$ for any such $i$, $j$, $e_i$, and $e_j$. Hence, \begin{align*}
\varrho(n) & = \varrho(\prod_{i=1}^{k}{p^{e_i}})\\
 & = \prod_{i=1}^{k}{\varrho( p^{e_i})}\\
 & = \prod_{i=1}^{k}{p^{\lfloor\frac{e_i}{2}\rfloor}}
\end{align*}
were the first equality follows by substitution of $n$ in its conical form, the second equality holds by $\S$ \ref{Lemma 3.2} because each factor is relatively prime, and the third holds by lemma $\S$ \ref{Lemma 3.1}.
\end{proof}
\end{Prop}
\section{Proper Congruent Numbers and the Altered Fermat Equation $ax^4+by^4 = z^2$}
\label{S4}
Fermat was the first to prove that the equation $x^4+y^4 = z^2$ has no-non trivial solution $(x,y,z)\in \mathbb{Z}^3$. In this section, I illustrate how \textit{proper} congruent numbers are related to the diophantine solutions of the equation $ax^4+by^4 = z^2$ for suitable integers $a$ and $b$. Specifically, I show that a whole integer $n\in \mathbb{N}$ being congruent relies heavily on the solutions of $ax^4+by^4=z^2$ were $\gcd(a,b)=1$, and $\mu(a)=0=\mu(b)$ or simply just one of them, that is $\mu(a)=0$ or $\mu(b)=0$. Furthermore, it is shown that $ax^4+by^4 = z^2$ has a solution if and only if $a+b=z^2$ which occurs if and only if there is a solution $c\in \mathbb{N}$ so that $\frac{1}{a}+\frac{1}{b}=\frac{1}{c}$ were $\gcd(a,b,c)$.
\begin{Prop}\label{Prop 4.1}
Each proper congruent number gives rise to a solution of the equation $ax^4+by^4=z^2$ were $\mu(a)=0=\mu(b)$ and $\gcd(a,b)=1=\gcd(c,d)$. \begin{proof} First assume that $n$ is a proper congruent number. By definition there exist $(\frac{a}{b},\frac{c}{d},\frac{e}{f})$ so that $n=\frac{1}{2}(\frac{a}{b})(\frac{c}{d})$ were as $(da)^2+(bc)^2 = e^2$ with the condition that $d\vert a$,$b\vert c$ and $d\neq 1\neq b$. This being the case we may write $a = dk$ and $c = bl$ for some $(k,l)$ so that $(da)^2+(bc)^2 = e^2$  becomes $(d^2k)^2+(b^2l)^2 = e^2$. This is a solution to $ax^4+by^4=z^2$ were $\mu(a)=0=\mu(b)$. By Proposition $\S$ \ref{Prop 2.3} we now that $\gcd(a,b)=1=\gcd(c,d)$. \end{proof}
\end{Prop}
\begin{Prop}\label{Prop 4.2}
Let $2n\in \mathbb{N}$ were $n>1$. If there is a divisor $d$ of $2n$ so that $(dx^2)^2+((\frac{2n}{d}y^2))^2=z^2$ for some $d\vert 2n$, then $n$ is congruent.\begin{proof}
Let $2n\in \mathbb{N}$ were $n>1$ and assume there is a divisor $d$ of $2n$ so that $(dx^2)^2+((\frac{2n}{d}y^2))^2=z^2$ for some $d\vert 2n$ and $x,y\in \mathbb{Z}^+$. Rearrange this to look like $(x(dx))^2+(y(\frac{n}{d}y))^2=z^2$. By Theorem $\S$ \ref{Thm 2.2} it follows that the number $\frac{1}{2}(\frac{dx}{y})(\frac{\frac{2n}{d}y}{x})$ is congruent. This is precisely $n$.
\end{proof}
\end{Prop}
\begin{Cor}\label{Corollary 4.1}
Let $p$ be prime. Then $p$ is congruent if and only if one of the equations $$x^4+4p^2y^4=z^2, ~4x^4+p^2y^4=z^2$$ has a solution for some $(x,y,z)\in \mathbb{N}^3$.\begin{proof}
Suppose that $p$ is prime. In view of proposition $\S$ \ref{Prop 4.2} if there is a divisor $d$ of $2p$ so that $(dx^2)^2+((\frac{2p}{d}y^2))^2=z^2$ for some $d\vert 2p$, then $p$ is congruent. Since $p$ is prime the only divisors are $d=1,2$ or $d=p$. Each case results in the three equations respectfully $$x^4+4p^2y^4=z^2, ~4x^4+p^2y^4=z^2,~p^2x^4+4y^4=z^2$$ Two of the equations are the same so this reduces to the equations in the hypothesis. Conversely, suppose one of these equations are satisfied. For any equation there is a $d\in \{1,2,p\}$ for which we may rewrite this as $(x(dx))^2+(y(\frac{2p}{d}y))^2=z^2$. Note that $x\vert dx$ and $y\vert (\frac{2p}{d}y)$ so by Theorem $\S$ \ref{Thm 2.2} the number $\frac{1}{2}\left(\frac{dx}{y}\right)\left(\frac{(\frac{2p}{d}y)}{x}\right)$ is congruent. This is precisely $p$.
\end{proof}
\end{Cor}
\begin{Thm}\label{Theorem 4.1}
Let $a,b\in \mathbb{N}$ were $a,b>1$ and $\gcd(a,b)=1$. Then the equation $ax^4+by^4=z^2$ has a solution $(x,y,z)\in \mathbb{N}^3$ if and only if the equation $a+b=z^2$ has a solution $z\in \mathbb{Z}$.\begin{proof} Let $a,b>1$ and assume the equation $ax^4+by^4=z^2$ has a solution $(x,y,z)\in \mathbb{N}^3$. Then we have the following implications 
\begin{align*}
ax^4+by^4&=z^2\\
& \Downarrow\\
\frac{a}{b}x^4+y^4 & = \frac{1}{b}z^2\\
& \Downarrow\\
\frac{a}{b}x^4 & = \frac{1}{b}z^2-y^4\\
& \Downarrow\\
\frac{a}{b} & = \frac{\frac{1}{b}z^2-y^4}{x^4}\\
& \Downarrow\\
\frac{a}{b} & = \frac{z^2-by^4}{bx^4}\\
\end{align*}
Now $b\vert by^4$ but $b\not \vert z^2$ since $\gcd(b,z)=1=\gcd(x,z)=\gcd(x,y)$. This means that $\gcd(z^2-by^4,bx^4)=1=\gcd(a,b)$. In other words, the fractions we found to be equal are actually in lowest terms. This implies that $b=bx^4$ so then $x=\pm 1$. By analogy we have the implication
\begin{align*}
ax^4+by^4&=z^2\\
& \Downarrow\\
x^4+\frac{b}{a}y^4 & = \frac{1}{a}z^2\\
& \Downarrow\\
\frac{b}{a} y^4 & = \frac{1}{b}z^2-x^4\\
& \Downarrow\\
\frac{b}{a} & = \frac{\frac{1}{b}z^2-x^4}{y^4}\\
& \Downarrow\\
\frac{b}{a} & = \frac{z^2-bx^4}{by^4}\\
\end{align*}
Again we see $a\vert ay^4$ but $a\not \vert z^2$ since $\gcd(a,z)=1=\gcd(y,z)=\gcd(x,y)$. This means that $\gcd(z^2-bx^4,ay^4)=1=\gcd(a,b)$. In other words, the fractions we found to be equal are actually in lowest terms. This implies that $b=bx^4$ so then $x=\pm 1$. Hence, we have found $x=1=y$. So then $ax^4+by^4=z^2$ implies that $a+b=z^2$.\\
\indent Secondly suppose that $a+b=z^2$ for some $z^2$. Then note that this is a solution to $ax^4+by^4=z^2$ were $x=\pm 1=y$. This completes the proof.
\end{proof}
\end{Thm}
\begin{Thm}\label{Theorem 4.2}
If there exist $a,b,c\in \mathbb{Z}$ so that $\frac{1}{a}+\frac{1}{b}=\frac{1}{c}$ were as $\gcd(a,b,c)=1$, then the the equation $a+b=z^2$ has a solution $z\in \mathbb{Z}$. \begin{proof}
Assume there exist a $a,b,c\in \mathbb{N}$ so that $\frac{1}{a}+\frac{1}{b}=\frac{1}{c}$ were as $\gcd(a,b,c)=1$. Multiple through by $abc$ to obtain $bc+ac=ab$. So then $c(b+a)=ab$. Let $\gcd(a,b)=d$ and set $a=a^{\prime}d$ and $b=b^{\prime}d$. Thus we have $\frac{1}{d}(c(b+a))=\frac{1}{d}ab$. This is equivalent to $c(b^{\prime}+a^{\prime})=b^{\prime}a^{\prime}d$. Since $a^{\prime}$, and $b^{\prime}$ are relatively prime so it must be that $(b^{\prime}+a^{\prime})\vert d$. But $d\vert (b^{\prime}+a^{\prime})$ so this means that $(b^{\prime}+a^{\prime}) = d$. Hence, $(a+b) = (a^{\prime}d+b^{\prime}d)=d(a^{\prime}+b^{\prime})=d^2$. 
\end{proof}
\end{Thm}
\begin{Prop}\label{Proposition 4.3}
Each proper congruent number gives rise to a solution of the equation $ax^4+y^4=z^2$ for some $a\in \mathbb{N}$. \begin{proof} By Proposition $\S$ \ref{Prop 4.1} each proper congruent number corresponds to an integer solution to the equation $ax^4+by^4=z^2$ for some relatively prime numbers $a$ and $b$. By Theorem $\S$ \ref{Theorem 4.2}, since $x$ and $y$ cannot be $1$ (this follows since $n$ is proper) it follows that $a=1$ or $b=1$, but they cannot both be $1$. In any either case were $a=1$ or $b=1$ we may rearrange the equation to look like $ax^4+y^4=z^2$ or $x^4+by^4=z^2$. Since addition is commutative and the fact that $x$ and $y$ are variables these equations are the same. Thus each proper congruent number corresponds to an integral solution to the equation $ax^4+y^4=z^2$. \end{proof}
\end{Prop}
\section{Semi-Proper Congruent Numbers}\label{S5}
\begin{Def}\label{Defintion 5.1}
Let $n$ be a congruent number which corresponds to triangle $\triangle ABC$ with rational sides $(\frac{a}{b},\frac{c}{d},\frac{e}{f})$. If $b=1$ or $d=1$ but not both, then we say that $n$ is a \textit{Semi-proper congruent number}.
\end{Def}
\begin{Prop}\label{Proposition 5.1}
If $n$ is a semi proper congruent number with $\triangle ABC$ constructed by the triple $(a,\frac{c}{d},\frac{e}{f})$, then $f=d$ and $e^2=c^2+ad^2$. \begin{proof}
Suppose $n$ is a semi proper congruent number with $\triangle ABC$ constructed by the triple $(a,\frac{c}{d},\frac{e}{f})$. Since $n$ is the area of the $\triangle ABC$ we have the identity $n = \frac{1}{2}(a)(\frac{c}{d})$ with $\gcd(c,d)=1$. Thus $d\not \vert c$ but this means that then $d\vert a$. Write $a=dk$ for some $k\in \mathbb{N}$. Thus 
\begin{align*}
(\frac{c}{d})^2+ a^2 &= (\frac{e}{f})^2\\
& \Downarrow \\
 (\frac{c}{d})^2 + (dk)^2 & = (\frac{e}{f})^2\\
& \Downarrow \\
 c^2+d^4k^2 & = d^2(\frac{e}{f})^2
\end{align*}
Recall that $\gcd(e,f)=1$ by assumption. Since $ c^2+d^4k^2 = d^2(\frac{e}{f})^2$ and the left side of the equation is a positive whole number it must follows that $f\vert d$. Since $f\vert d$ there exist a $l$ such that $d=fl$. Hence, we have the following implication:
\begin{align*}
(\frac{c}{d})^2+ a^2 &= (\frac{e}{f})^2\\
& \Downarrow\\
 (\frac{c}{fl})^2+ a^2 &= (\frac{e}{f})^2\\
 & \Downarrow\\
(\frac{c^2}{(fl)^2}+ \frac{a^2(fl)^2}{(fl)^2} &= (\frac{e}{f})^2\\
& \Downarrow\\
(\frac{c^2+a^2(fl)^2}{(fl)^2}) &= (\frac{e}{f})^2\\
\end{align*}
By assumption we know that $\gcd(e,f)=1$. Therefore, the right side of the expression is in lowest terms. Since $d=fl\not \vert c$  it follows that the left side of the expression is in lowest term. Two quotients in lowest terms have equal numerator and denominator. Hence, $f^2=f^2l^2$ which implies $l=1$. Since $f=dl$ and $l=1$ we find $f=d$ and $e^2=c^2+ad^2$ as desired.
\end{proof}
\end{Prop}
\begin{Thm}\label{Theorem 5.1}
If $a,c,e,d$ are integers that satisfy the following conditions:\\
\noindent  (i) $e^2 = c^2+(da)^2$\\
(ii) $\gcd(c,d) = 1$ and $d\vert a$\\
(iii) $\{a,c\}$ contains an even integer,  and this even number has more factors of 2 than $2d$ has.\\
Then the number $\frac{1}{2}(a)(\frac{c}{d})$ is a Semi-proper congruent number. \begin{proof}
Let conditions $(i)-(iii)$ hold valid. Set $f=d$, then we have $(\frac{e}{f})^2 = (\frac{c}{d})^2+a^2$. This is a Pythagorean triple which represents a triangle with area $\frac{1}{2}(a)(\frac{c}{d})$. Thus this integer is a Semi-proper congruent number.
\end{proof}
\end{Thm}
\begin{Prop}\label{Proposition 5.2}
Let $a>1$ and suppose $x^2+a^2y^4 =z^2$ for some $x,y\in \mathbb{N}$. Then the number $n=\frac{1}{2}ax$ is a Semi-Proper congruent number. \begin{proof}
Suppose that $a>1$ and $x^2+a^2y^4 =z^2$ for some $x,y\in \mathbb{N}$. Rewrite this as $x^2+(y(ay))^2 =z^2$. Since this is a Pythagorean triple by proposition $\S$ \ref{Prop 2.3} it follows that  $\gcd(x,y(ay))=1$. Since $y\vert ay$ and at least one of $\{x,ay\}$ is even conditions $(i)-(iii)$ of Theorem $\S$ \ref{Theorem 5.1} are satisfied. As a consequence the number $\frac{1}{2}(ay)(\frac{x}{y})=\frac{1}{2}ax$ is semi-congruent.
\end{proof}
\end{Prop}
\section{Improper Congruent Numbers}\label{S6}
This section is short because there is not to much to prove or show. Improper congruent numbers are simple congruent numbers because they represent the are of a right triangle whose sides are specifically whole integers and not rational numbers. I call these simple because Euclid has already described these type of numbers for us.
 \begin{Def}\label{Def 6.1}
A \textit{Improper Congruent} number is a number $n\in \mathbb{N}$ so that $n$ is the area of some right triangle $\triangle ABC$ with integer sides $(a,b,c)$ That is all sides lengths of the $\triangle ABC$ are strictly in $\mathbb{Z}^{+}$ and not in $\mathbb{Q}$.
\end{Def}
\begin{Thm} 
\label{Theorem 6.1} 
Let $k\in \mathbb{N}$. If $k$ is a improper congruent number then it can be written as $mn(m-n)(m+n)$ for some $m,n \in \mathbb{Z}^+$ such that $\gcd(m,n)=1$, $m>n$ and $m \not \equiv n(\bmod{2})$. Moreover, $mn(m-n)(m+n)$ is also some improper congruent number. \begin{proof}
By Theorem $\S$ \ref{Thm 2.1} it follows that all primitive Pythagorean triples can be written as 
\begin{align*}
a & = m^2-n^2\\
b& = 2nm\\
c &= m^2+n^2
\end{align*}
were $\gcd(n,m)=1$ and $m>n$ and $m\not \equiv n(\bmod{2})$. The area of each such triangle is thus $k = \frac{1}{2}2mn(m-n)(m+n) = mn(m-n)(m+n)$
\end{proof}
\end{Thm}
\section{Characterizing Congruent Numbers Using the Mobius function}\label{S7}
In this section I show how to generate all the types of congruent number using Pythagorean triples. We have  already done two examples of this in Ex. $\S$ \ref{Ex 2.1}, and Ex. $\S$ \ref{Ex 2.1}. I will show that all congruent numbers can be written as $$\frac{nm(m-n)(m+n)}{\left(\tau\sigma\right)^2}$$ were as $\tau\vert \varrho(2nm)$ and $\sigma\vert \varrho((m-n)(m + n))$ for some $m,n$. I will show that Improper congruent numbers correspond to this representation when the values $d_1 = 1 = d_2$, Semi-proper congruent numbers correspond to only one value $\tau = 1$ or else $\sigma = 1$. Moreover, I demonstrate that each proper congruent number corresponds to values $\tau\neq 1 \neq  \sigma$. In order to accomplish this I use the preliminaries that were established throughout this work. As a consequence we develop some useful results about congruent numbers. Specifically, these results show that some of that the divisors $d$ of $\varrho(mn(m-n)(m+n)$ and every congruent number corresponds to one of these divisors. 
\begin{Thm}\label{Theorem 7.1}
Every Improper-Congruent number can be written as $m^3n-n^3m$ whereas $m$ is not square free, $n$ is not square free,or possibly when both $m$ and $n$ are not square free, and were as $\gcd(m,n) = 1$, $m > n$, and $m\neq n (\bmod{2})$.
\begin{proof}
Suppose that $\alpha$ is an Improper-Congruent number. By definition $\alpha$ is the area of a right triangle  $\triangle ABC$ with all sides $\{a,b,c\}$ all in $\mathbb{Z}^+$. By Theorem $\S$ \ref{Thm 2.1} there exist $m$ and $n$ such that \begin{align*}
a & = m^2-n^2\\
b& = 2nm\\
c &= m^2+n^2
\end{align*}
The area of this triangle $\triangle ABC$ is $\frac{1}{2}\left( m^2-n^2\right)\left(2nm\right) = m^3n-n^3m$.
\end{proof}
\end{Thm} 
\begin{Thm}\label{Theorem 7.2}
Every Semi-proper congruent number can be written as  $$\frac{mn(m-n)(m+n)}{\sigma^2}$$ were as $\sigma\vert \varrho(nm(m-n)(m+n))$,and $\gcd(m,n) = 1$, $m > n$, and $m\neq (\bmod{2})$ with $\sigma\vert nm$ or $\sigma\vert (m-n)(m+n) $.
\begin{proof}
The proof follows by an argument of decent. Suppose that $\alpha$ is an Semi-Proper congruent number corresponding to the triple $\{a,\frac{c}{d},\frac{e}{f}\}$ with and $\mu(\alpha) = 1$. Then by proposition $\S$ \ref{Proposition 5.1}  it follows that $e^2=c^2+a^2d^2$. Note that this is a solution to the equation $x^2+a^2y^2=z^2$ with $x=a$, $y=d$, and $z=e$. Now suppose that $\mu(y) =0$. Then there exist a $d\vert y$ such that $d^2\vert y$. Hence, write $y=d^2k_1$ and notice that the equation $x^2+a^2y^2=z^2$ transforms into the equation $x^2+(ak_1)^2d^4=z^2$. Note that this is a solution to the equation $x^2+a^2 y^4=z^2$. Thus by direct application of proposition $\S$ \ref{Proposition 5.2} it follows that the number $\frac{1}{2}xak_1$ is a Semi-Congruent number. Now we apply the same argument again. If $\mu(d)=0$ then there exist a $\tau\vert d$ such that $\tau^2\vert d$. Therefore we may write $d = \tau^2 k_2$ for some $k_2\in \mathbb{Z}$. Then by direct substitution the equation $x^2+(ak_1)^2d^4=z^2$ transforms into the equation $$x^2+(ak_1k_2)^2\tau^4=z^2$$ and again by proposition $\S$ \ref{Proposition 5.2} it follows that $\frac{1}{2}xak_1k_2$ is a Semi-Congruent number. We take note that if $k_1\neq 1\neq k_2$ then $\frac{1}{2}xak_1<\frac{1}{2}xak_1k_2$, which implies that these are two distinct Semi-Congruent numbers. By repeating this argument until we reach a point were the equation $x^2+a^2 y^4=z^2$ has $\mu(y) = 0$ we find a set of Semi-Congruent numbers, which all correspond to this single Pythagorean triple. In complete generality each congruent number we find looks like $\frac{1}{2}x\tau$ for some $\tau\vert (ay)^2 $ with $\mu(t) = 0$. Hence, since $\tau$ is a divisor of the number $(ay)^2$ it must divide $ay$. By Euclid's theorem $\S$ \ref{Thm 2.1} we have $ay = (m-n)(m+n)$ or else $ay = 2nm$ for some $m$ and $n$ relatively prime and of opposite parity. Hence, $\tau$ must divide $(m-n)(m+n)$ or else $2nm$. But then $\tau$ can be written as $\frac{2mn}{\sigma^2}$ or $\frac{(m-n)(m+n)}{\sigma^2}$ for some $\sigma \vert (m-n)(m+n) $ or $\sigma\vert 2mn$. In either case the are of such a triangle is $$\frac{1}{2}x\tau = \frac{1}{2}\frac{2mn(m-n)(m+n))}{\sigma^2} = \frac{nm(m-n)(m+n)}{\sigma^2}$$ In the second case when $mu(y)=1$ we find $\sigma = 1$ yields the result.
Since, this holds for each divisor $\tau$ with $\mu(t) = 0$ the result of the hypothesis then follows.
\end{proof}
\end{Thm} 
\begin{Thm}\label{Theorem 7.3}
Every Proper congruent number can be written as  $$\frac{mn(m-n)(m+n)}{(\sigma \tau)^2}$$ were as $\sigma^2\vert \varrho(mn)$, $\tau\vert \varrho((m-n)(m+n))$ and $\gcd(m,n) = 1$, $m > n$, and $m\neq n(\bmod{2})$.
\begin{proof}
The proof of this theorem is similar to the proof of the last two theorems. The proof is summarized and not given rigorously to make sure the reader has paid close attention thus far. First suppose that $\alpha$ is a Proper-Congruent number. Then we want to show that it can be written in the form presented in the hypothesis. First since $\alpha$ is Proper it represents the area of some triangle corresponding to the triple $\{\frac{a}{b},\frac{c}{d},\frac{e}{f}\}$. So then $e^2 = (da)^2+(bc)^2$. Now this a Pythagorean triple which corresponds a triangle with areas $\frac{1}{2}\frac{a}{b}\frac{c}{d}$. The trick here is we can reformulate $d$, $a$, $b$, and $c$ by simply writing the factorization of them all. However, we notice that of we shift the prime factorization of $d$ into $a$, meaning we take some prime from $d$ and give them to $a$ and write $da = (\tau a)$ were $\frac{a}{d} = \tau$. Then we see that $e^2 = (\tau a)^2+(bc)^2 $. As a consequence of Theorem $\S$ \ref{Thm 2.1} it follows that the number $n = \frac{1}{2}\frac{a}{b}\frac{c}{\tau}$ is also a Proper-Congruent number. Clearly, the quantities $\frac{1}{2}\frac{a}{b}\frac{c}{d}$ and $\frac{1}{2}\frac{a}{b}\frac{c}{\tau}$ are not the same if $\tau \neq d$. Now we can repeat the argument until we exhaust all the non square free divisors of $a$, and then use the same argument for the quantity $b$, and both $a$ and $b$ simultaneously. By doing this we find the result, by finding the area of each each one of these triangles. These argument is summarized by noting that the number of congruent number that correspond to each Pythagorean triple is equal to the number of non square free divisors of for its components $a$ and $b$ corresponding to the triple $\{a,b,c\}$.
\end{proof}
\end{Thm} 
\begin{Thm}\label{Theorem 7.4}
Let $A$ be given as  $$A = \{n\in \mathbb{N}: \text{n is a proper congruent number }\}$$ and $B$ be given as $$B = \{(a,b,c) \in \mathbb{N}^{3}: a^2+b^2 = c^2, \mu(a)=0=\mu(b)\}$$ Let $(\alpha,\beta,\gamma)\in B$. and write $\alpha$, and $\beta$ in canonical form: $$\alpha = \prod_{p\vert a}^{}{p^{e_i}},\indent \beta = \prod_{p\vert b}^{}{p^{e_i^{\prime}}} $$ Then the number of congruent numbers that correspond to the triple $(\alpha,\beta,\gamma)$ is given by the quantity $$\sum_{i=1}^{\max(k,k^{\prime})}{\lfloor\frac{e_i+e_i^{\prime}}{2}\rfloor}$$ were $k$ and $k^{\prime}$ denotes the number of prime factors of $\alpha$ and $\beta$ respectfully. \begin{proof}
This proof is left to the reader as an exercise. (Hint! Use Table 3.1 to help understand the hypothesis, then explain why each non-square divisor of $\alpha$ and $\beta$ corresponds to a unique congruent number and after add them both.)
\end{proof}
\end{Thm}


\newpage
\begin{thebibliography}{00}
\bigskip
%
\bibitem{alvaro}{\textcolor{blue}{PLozano-Robledo, Alvaro, and Alvaro Lozano-Robledo. Elliptic curves, modular forms, and their L-functions. Providence, RI: American Mathematical Society, 2011.}}
\bibitem{m184}{\textcolor{blue}{ Conrad, Keith. "The congruent number problem." The Harvard College Mathematics Review 2 (2008): 58-74.}}
\bibitem{RFEE}{\textcolor{blue}{Fitzpatrick, Richard. Euclid's Elements In Greek. Lulu.com, 2006}}
%
\bibitem{alvero}{\textcolor{blue}{ Lozano-Robledo, Alvaro, and Alvaro Lozano-Robledo. Elliptic curves, modular forms, and their L-functions. Providence, RI: American Mathematical Society, 2011.}}
%
\bibitem{GI}{\textcolor{blue}{Koblitz, Neal I. Introduction to elliptic curves and modular forms. Vol. 97. Springer Science \& Business Media, 2012.}}
%
\bibitem{Wiles}{\textcolor{blue}{Wiles, A. "The Birch and Swinnerton-Dyer Conjecture, Clay Mathematics Institute Web Site."}}
%
\bibitem{Tunnel}{\textcolor{blue}{J. Tunnell: A Classical Diophantine Problem and Modular Forms of Weight 3/2, Invent. Math.72(1983),323–334}}
%
\bibitem{Conrad}{\textcolor{blue}{Conrad, Keith. "The congruent number problem." The Harvard College Mathematics Review 2 (2008): 58-74.}}
%
\bibitem{Arab}{\textcolor{blue}{Tian, Ye. "Congruent numbers with many prime factors." Proceedings of the National Academy of Sciences 109.52 (2012): 21256-21258.}}
\end{thebibliography}
\end{document}